\documentstyle[12pt]{article}
\makeatletter
\renewcommand\thesection{\@Roman\c@section}
\renewcommand\thesubsection{\thesection.\@arabic\c@subsection}
\makeatother

\begin{document}


\title{\hspace{11cm}{\small\bf IMPNWU-980108}\\
\vspace{2cm}
{\Large \bf The dynamical twisting and nondynamical r-matrix structure of  
elliptic Ruijsenaars-Schneider model}}

\author{
 Bo-yu Hou$^{a,b}$ and Wen-Li Yang$^{a,b}$
\thanks{e-mail :wlyang@phy.nwu.edu.cn}
\thanks{Fax    :86-29-8302331}
\\
\bigskip\\
$^{a}$ CCAST ( World Laboratory ), P.O.Box 8730 ,Beijing 100080, China\\
$^{b}$ Institute of Modern Physics, Northwest
University, Xian 710069, China   
\thanks{Mailing address}}

\date{}
\maketitle
\begin{abstract}
From the dynamical twisting of the classical r-matrix, we obtain a new Lax
operator for the elliptic Ruijsenaars-Schneider model (cf. Ruijsenaars').
The
corresponding r-matrix is shown to be the  classical $Z_n$-symmetric
elliptic r-matrix, which is the same as that obtained in the study of the
nonrelativistic version---the $A_{n-1}$ Calogero-Moser model.

\vspace{0.6truecm}
\end{abstract}


\def\beq{\begin{equation}}
\def\eeq{\end{equation}}
\def\bea{\begin{eqnarray}}
\def\eea{\end{eqnarray}}  
\def\ba{\begin{array}}  
\def\ea{\end{array}}
\def\no{\nonumber}
\def\lt{\left}
\def\rt{\right}
\newcommand{\bq}{\begin{quote}}
\newcommand{\eq}{\end{quote}}
\newtheorem{Theorem}{Theorem}
\newtheorem{Definition}{Definition}
\newtheorem{Proposition}{Proposition}
\newtheorem{Lemma}{Lemma}
\newtheorem{Corollary}{Corollary}
\newcommand{\proof}[1]{{\bf Proof. }
        #1\begin{flushright}$\Box$\end{flushright}}

\newcommand{\sect}[1]{\setcounter{equation}{0}\section{#1}}
\renewcommand{\theequation}{\thesection.\arabic{equation}}

\section{Introduction}

Following the successes of the Calogero-Moser (CM) models
\cite{Cal75,Mos75}, a relativistic
generalization of the CM models--- the so-called Ruijsenaars-Schneider
(RS) models have been proposed \cite{Rui86}, which  the intergrability
has been conserved. The RS model describes a completely integrable system
of n one-dimensional interacting
relativistic particles. Its importance lies in the fact that it is related
to the dynamics of solitons in some integrable relativistic field
theories\cite{Bab93,Bra98} and its discrete-time
version has been connected with the Bethe anstaz equation of the solvable
lattice statistical model \cite{Nij96}. Recent development was
shown that it can be obtained by a Hamiltonian reducation of the cotangent
bundle of some Lie group \cite{Aru96}, and can also be considered as
the gauged WZW
theory \cite{Gor94}. The study of RS model would
play an universal role in study of completely integrable multi-particle
systems. Among all type RS models, the elliptic RS model is the most
general one and the other type such as the rational, hyperbolic and
trigonomettic type is just the various degenerations of the elliptic one. 
In this paper, we shall study the elliptic $A_{n-1}$ type RS model with 
generic n ( $n>2$).

The Lax representation and its corresponding r-matrix structure for
rational, hyperbolic and trigonometric $A_{n-1}$ type RS models were 
constructed by Avan et al \cite{Ava93}. The Lax representation for the
elliptic RS models was constructed by Ruijsenaars \cite{Rui88}
, and the corresponding r-matrix structure was given by Nijhoff et al 
\cite{Nij961} and Suris \cite{Sur96}. It turns out  that the
r-matrix structure of the 
RS model is given in terms of a quadratic Poisson-Lie bracket
 with dynamical r-matrices (i.e. the r-matrix  depends upon the dynamical
variables). Particularly, in contrast with
the dynamical Yang-Baxter equation of the r-matrix structure of  the CM
model, the generalized Yang-Baxter relations for the
quadratic Poisson-Lie bracket with a dynamical r-matrix is still an open
problem \cite{Nij961}. Since the Poisson bracket of the Lax operator is no
longer closed, the quantum version of such classical L-operator has not
been able to be constructed.

It is well-known that
the Lax representation for a completely integrable models is not unique.
It has been recongized \cite{Bab89,Hou97} that the r-matrix of a
model can be changed
drastically by the choice of Lax representation. 
In our former work \cite{Hou97}, we succeeded in constructing a new Lax
operator
(cf. Krichever's \cite{Kri80}) for the elliptic $A_{n-1}$ CM model
and showing that the corresponding r-matrix is a nondynamical one, which is
the classical $Z_n$-symmetric elliptic r-matrix \cite{Hou89,Hou97}.
Very recently, we
found a ``good" Lax operator for the elliptic RS model with a very
special case n=2 \cite{Hou99}. In present paper, extending  our former
work in \cite{Hou99}, we construct a ``good"  Lax operator (in such a
sense that it has a nondynamical r-matrix structure)  for the elliptic RS 
model with a general $n$ $(n>2)$.

The paper is organized as follows. In section 2, we construct the dynamical 
twisting relations of the classical r-matrix for the quadratic Poisson-Lie 
bracket. The condition that the `` good" Lax representation
could exist is found. In section 3, some briefly reviews of Nijhoff et al
's work  on dynamical r-matrix of the elliptic RS model was given. 
In section 4, we
construct the ``good"  Lax  representation for elliptic RS model with generic 
n, and obtain the corresponding nondynamical r-matrix structure. The
quantum version  L-operator of the Lax operator  is 
constructed in section 5. Finally, we give
summary and discussions. Appendix contains some detailed
calculations.

\section{The dynamical twisting of classical r-matrix}
In this section we will give  some general theories of
the completely integrable finite particles systems.

A Lax pair (L,M) consists of two functions on the phase space of the system
with values in some Lie algebra {\bf $g$}, such that the evolution equations
may be written in the following form
\begin{eqnarray}
\frac{dL}{dt}=[L,M],
\end{eqnarray}
where $[,]$ denotes the bracket in the Lie algebra {\bf $g$}. The
interest in the existence of such a pair lies in the fact that it allows
for an easy construction of conserved quantities (integrals of motion). 
It follows that the adjoint-invariant quantities $trL^{l}\ \ (l=1,...,n)$
are the integrals
of the motion. In order to implement Liouville theorem onto this set of
possible action variables we need them to be Poisson-commuting. As shown
in \cite{Bab89},  the commutativity of the integrals $trL^{l}$ of the
Lax
operator can be deuced from  that the fundamental Poisson bracket
$\{L_{1}(u),L_{2}(v)\}$ could be represented in the linear commutator form
\begin{eqnarray}
\{L_1(u),L_2(v)\}=[r_{12}(u,v),L_1(u)]-[r_{21}(v,u),L_{2}(v)],\label{Linear}
\end{eqnarray}
or quadratic  form \cite{Fre91}
\begin{eqnarray}
\{L_1(u),L_2(v)\}&=&L_1(u)L_2(v)r^{-}_{12}(u,v)-r^{+}_{21}(v,u)L_1(u)L_{2}(v)
\no\\
&&~~+L_1(u)s_{12}^{+}(u,v)L_2(v)-L_2(v)s_{12}^{-}(u,v)L_1(u),\label{Qua}
\end{eqnarray}
where we have used the notation
\begin{eqnarray*}
L_1\equiv L\otimes 1\ \ ,\ \ L_2\equiv 1\otimes L\ \ ,\ \
a_{21}=Pa_{12}P, 
\end{eqnarray*}
and $P$ is the permutation operator such that
$Px\otimes y=y\otimes x$.

The dynamical twisting of the linear Poisson-Lie bracket (\ref{Linear}) 
was studied in the \cite{Hou97} (we refer therein) and also studied by
Babelon et al \cite{Bab89}. We are  to
inverstigate the general dynamical twisting of the quadratic Poisson-Lie
bracket (\ref{Qua}). 

In order to define a consistent Poisson bracket, one should impose
some constraints on the r-matrices. The skew-symmetry of Poisson bracket
require that 
\begin{eqnarray}
& &r^{\pm}_{21}(v,u)=-r^{\pm}_{12}(u,v)\ \ ,\ \ s^{+}_{21}(v,u)
=s^{-}_{12}(u,v),\label{Skew1}\\
&
&r^{+}_{12}(u,v)-s^{+}_{12}(u,v)=r^{-}_{12}(u,v)-s^{-}_{12}(u,v).
\label{Skew2}
\end{eqnarray}
As for  the numerical r-matrices $r^{\pm}(u,v),s^{\pm}(u,v)$ case, 
some constraints condition (sufficient conidtion) imposed on the r-matrices
to make Jacobi identity satisfied, was given by Freidel et al
\cite{Fre91}. 
However, generally speaking, 
the Jacobi identity for the dynamical  
r-matrices $r^{\pm}(u,v),s^{\pm}(u,v)$  would take a very complicated
form. 

It should be remarked that such a classification (from  dynamical and
nondynamical r-matrix structure) is by no means unique, which
drastically depend on the Lax representation which  one choose for a  
system. Therefore, there is no one-to-one correspondence between a given
dynamical system and a defined r-matrix. The  same dynamical system may
have several Lax representations and several r-matrix \cite{Hou97}.
The different
Lax representation of a system is conjugated each other. Namely,  if
($\stackrel{\sim}{L},\stackrel{\sim}{M}$) is one of other Lax pair of the
same dynamical system conjugated with the old one $(L,M)$, it means that 
\begin{eqnarray}
\stackrel{\sim}{L}(u)=g(u)L(u)g^{-1}(u),~~~~
\stackrel{\sim}{M}(u)=g(u)M(u)g^{-1}(u)
-(\frac{d}{dt}g(u))g^{-1}(u),\label{Con}
\end{eqnarray}
where $g(u)\in G$ whose Lie algebra is {\bf $g$}. Then, we have

\noindent {\it {\bf Proposition 1.} The Lax pair
($\stackrel{\sim}{L},\stackrel{\sim}{M})$ has the following r-matrix
structure
\bea
\{\stackrel{\sim}{L}_1(u),\stackrel{\sim}{L}_2(v)\}&=&
\stackrel{\sim}{L}_1(u)\stackrel{\sim}{L}_2(v)\stackrel{\sim}{r}^{-}_{12}
(u,v)-\stackrel{\sim}{r}^{+}_{21}(v,u)
\stackrel{\sim}{L}_1(u)\stackrel{\sim}{L}_2(v)\no\\
&&
+\stackrel{\sim}{L}_1(u)\stackrel{\sim}{s}^{+}_{12}(u,v)
\stackrel{\sim}{L}_2(v)-
\stackrel{\sim}{L}_2(v)\stackrel{\sim}{s}^{-}_{12}(u,v)
\stackrel{\sim}{L}_1(u),\label{3a}
\eea
\noindent where
\bea
\stackrel{\sim}{r}^{-}_{12}(u,v)&=&g_1(u)g_2(v)r^{-}_{12}(u,v)
g^{-1}_1(u)g^{-1}_2(v)-\stackrel{\sim}{\Delta}_{12}(u,v)
+\stackrel{\sim}{\Delta}_{21}(v,u),\nonumber\\
\stackrel{\sim}{r}^{+}_{12}(u,v)&=&g_1(u)g_2(v)r^{+}_{12}(u,v)
g^{-1}_1(u)g^{-1}_2(v)-\stackrel{\sim}{\Delta}^{(1)}_{12}(u,v)
+\stackrel{\sim}{\Delta}^{(1)}_{21}(v,u),\nonumber\\
\stackrel{\sim}{s}^{+}_{12}(u,v)&=&g_1(u)g_2(v)s^{+}_{12}(u,v)
g^{-1}_1(u)g^{-1}_2(v)-\stackrel{\sim}{\Delta}_{21}(v,u)
-\stackrel{\sim}{\Delta}^{(1)}_{12}(u,v),\nonumber\\
\stackrel{\sim}{s}^{-}_{12}(u,v)&=&g_1(u)g_2(v)s^{-}_{12}(u,v)
g^{-1}_1(u)g^{-1}_2(v)-\stackrel{\sim}{\Delta}_{12}(u,v)
-\stackrel{\sim}{\Delta}^{(1)}_{21}(v,u),\nonumber\\
\stackrel{\sim}{\Delta}_{12}(u,v)&=&\stackrel{\sim}{L}^{-1}_{2}(v)
\Delta_{12}(u,v),~~~~~~~~~~~\stackrel{\sim}{\Delta}^{(1)}_{12}(u,v)=
\Delta_{12}(u,v)\stackrel{\sim}{L}^{-1}_{2}(v), \nonumber\\
\Delta_{12}(u,v)&=&
+\frac{1}{2}[\{g_1(u),g_2(v)\}g^{-1}_1(u)g^{-1}_2(v),
g_2(v)L_2(v)g^{-1}_2(v)]\no\\
& &\times 
g_{2}(v)\{g_{1}(u),L_2(v)\}g^{-1}_1(u)g^{-1}_2(v)
\nonumber
\eea
and the properties of (\ref{Skew1}) and (\ref{Skew2})  are conserved
\begin{eqnarray}
& &\stackrel{\sim}{r}^{\pm}_{21}(v,u)=-\stackrel{\sim}{r}^{\pm}_{12}(u,v)
\ \ ,\ \ \stackrel{\sim}{s}^{+}_{21}(v,u)
=\stackrel{\sim}{s}^{-}_{12}(u,v),\no\\
& &\stackrel{\sim}{r}^{+}_{12}(u,v)-\stackrel{\sim}{s}^{+}_{12}(u,v)=
\stackrel{\sim}{r}^{-}_{12}(u,v)-\stackrel{\sim}{s}^{-}_{12}(u,v).\no
\end{eqnarray} }
\vspace{0.6cm}

\noindent {\it \bf Proof:} The proof  is direct substituting 
(\ref{Con}) into the
fundamental Poisson bracket (\ref{Qua}) and use the following identity
\begin{eqnarray*}
\left[ [a_{12},L_1],L_2\right] =\left[ [a_{12},L_2],L_1\right].
\end{eqnarray*}
where $a_{12}$ is any matrix on {\bf $g\otimes g$}. 
\hspace{3cm} ${\bf\large\Box}$

It can be seen that: {\bf I.} The Lax operator $L(u)$ is transfered
as a similarity transformation from the different Lax representation; 
{\bf II.} The corresponding $M$ is undergone the usual gauge
transformation; {\bf III.} The r-matrices are transfered as some
generalized gauge transformation, which can be considered as 
the generalized
classical version of the dynamically twisting relation of  the quantum
R-matrix \cite{Ava96}. Therefore, it is of great value to find a
$``$good" Lax
representation for a system if it exists,
in which the corresponding r-matrices are all  nondynamical ones and
$r^{+}_{12}(u,v)=r^{-}_{12}(u,v)\ \
,\ \ s^{\pm}_{12}(u,v)=0$. In this special case, the corresponding
Poisson-Lie bracket becomes the Sklyanin bracket and  the well-studied
theories \cite{Bel84,Fad87} can be directly applied in the system.

\noindent {\it {\bf Corollary  1. } For given Lax pair $(L,M)$ and
the corresponding r-matrices, if there exist  $g(u)$ satisfied 
\bea
& &g_1(u)g_2(v)s^{+}_{12}(u,v)g^{-1}_1(u)g^{-1}_2(v)
-\stackrel{\sim}{\Delta}_{21}(v,u)-\stackrel{\sim}{\Delta}^{(1)}_{12}(u,v)=0
,\nonumber\\
& &g_1(u)g_2(v)s^{-}_{12}(u,v)
g^{-1}_1(u)g^{-1}_2(v)-\stackrel{\sim}{\Delta}_{12}(u,v)
-\stackrel{\sim}{\Delta}^{(1)}_{21}(v,u)=0,\label{Con1}\\
& &
\partial_{q_i}h_{12}=\partial_{p_j}h_{12}=0,\label{Con2}
\eea
where 
\bea
& &h_{12}(u,v)=g_1(u)g_2(v)r^{-}_{12}(u,v)g^{-1}_1(u)g^{-1}_2(v)
-\stackrel{\sim}{\Delta}_{12}(u,v)+\stackrel{\sim}{\Delta}_{21}(v,u)
\nonumber\\
& &\ \ \ \ \equiv g_1(u)g_2(v)r^{+}_{12}(u,v)
g^{-1}_1(u)g^{-1}_2(v)-\stackrel{\sim}{\Delta}^{(1)}_{12}(u,v)
+\stackrel{\sim}{\Delta}^{(1)}_{21}(v,u),
\end{eqnarray}
the nondynamical  Lax representation with Sklyanin Poisson-Lie
Bracket of the system would exist. }

\vspace{0.6cm}

\noindent The main purpose of this paper is  to find a
``good" Lax representation for the elliptic RS model with generic
$n$ $(n>2)$.

\section{Review of  elliptic RS model}
We first define  
some elliptic functions
\begin{eqnarray}
& &\theta^{(j)}(u)=
\theta\left[\begin{array}{c}\frac{1}{2}-\frac{j}{n}\\ 
\frac{1}{2}\end{array}\right](u,n\tau),~~~
\sigma(u)=\theta\left[\begin{array}{c}\frac{1}{2}\\
\frac{1}{2}\end{array}\right](u,\tau),\\
& &\theta\left[\begin{array}{c}a\\ b\end{array}\right](u,\tau)
=\sum_{m=-\infty}^{\infty}exp\{\sqrt{-1}\pi[(m+a)^{2}\tau +
2(m+a)(z+b)]\},
\nonumber\\
& &\theta'^{(j)}(u)=\partial_{u}\{\theta^{(j)}(u)\}\ \ ,\ \
\sigma'(u)=\partial_{u}\{\sigma(u)\},~~~
\xi(u)=\partial_{u}\{ln\sigma(u)\},
\end{eqnarray}
\noindent where $\tau$ is a complex number with $Im(\tau)>0$.

The Ruijsenaars-Schneider model is the system of n one-dimensional 
relativistic  particles interacting by the two-body potential.
In terms of the canonical variables $p_i,q_i\ \ (i=1,...n)$ enjoying in
the canonical Poisson bracket
\begin{eqnarray*}
\{p_i,p_j\}=0,\ \ \{q_i,q_j\}=0,\ \ \{q_i,p_j\}=\delta_{ij},
\end{eqnarray*}
the Hamiltonian of the system is expressed as \cite{Rui88}
\begin{eqnarray}
H=mc^2\sum_{j=1}^{n}coshp_j\prod_{k\neq j}
\left\{\frac{\sigma(q_{jk}+\gamma)\sigma(q_{jk}-\gamma)}{\sigma^2(q_{jk})}
\right\}^{\frac{1}{2}},~~~
q_{jk}=q_j-q_k.\label{Ham}
\end{eqnarray}
Here, $m$ denotes the particle mass, $c$ denotes the speed of
light, $\gamma $ is the coupling constant. 
The Hamiltonian (\ref{Ham}) is known to be
completely integrable. The most effective way to show its
integrability is to construct the Lax representation for the system (namely,
 to find the classical Lax operator). One L-operator  for the elliptic RS
model was  given by Ruijsenaars \cite{Rui88}
 \begin{eqnarray}
 L_{R}(u)^{i}_{j}=\frac{e^{p_j}\sigma(\gamma+u+q_{ji})}{\sigma(\gamma+q_{ji})
 \sigma(u)}\prod_{k\neq j}^{n}
\left\{\frac{\sigma(q_{jk}+\gamma)\sigma(q_{jk}-\gamma)}
{\sigma^2(q_{jk})}\right\}^{\frac{1}{2}} \ \ ,\ \ i,j=1,...,n.\label{Ru}
\end{eqnarray}
Alternatively, we adopt another Lax operator $\stackrel{\sim}{L}_{R}$, 
which is similar to that of Nijhoff et al in \cite{Nij961}
\begin{eqnarray}
\stackrel{\sim}{L}_{R}(u)^{i}_{j}=\frac{e^{p_j}\sigma(\gamma+u+q_{ji})}
{\sigma(u)\sigma(\gamma+q_{ji}}
\prod_{k \neq j}\frac{\sigma(\gamma +q_{jk})}{\sigma(q_{jk})}.\label{Lr}
\end{eqnarray}
The relation of $\stackrel{\sim}{L}_{R}$ with the standard 
Ruijsenaars'  $L_{R}(u)$ can be obtained from a Poisson map (or a canonical 
transformation)
\begin{eqnarray}
q_i\longrightarrow q_i\ \ \ ,\ \ p_i\longrightarrow p_i+\frac{1}{2}
ln\prod_{k\neq
i}\frac{\sigma(q_{ik}+\gamma)}{\sigma(q_{ik}-\gamma)}.\label{Map}
\end{eqnarray}

\noindent {\it {\bf Proposition 2.} The map defined in (\ref{Map}) is a
Poisson map.}

\noindent {\it\bf Proof:} The proposition 2 can be proven from considering
the symplectic two-form
\begin{eqnarray*}
& &\sum_{i}d(p_i+\frac{1}{2}
ln\prod_{k\neq i}\frac{\sigma(q_{ik}+\gamma)}{\sigma(q_{ik}-\gamma)})\wedge 
dq_i\\
& &\ \ \ =\sum_{i}dp_i\wedge dq_i -\frac{1}{2}\sum_{k\neq i}
(\frac{\sigma'(q_{ik}+\gamma)}{\sigma(q_{ik}+\gamma)}
-\frac{\sigma'(q_{ik}-\gamma)}{\sigma(q_{ik}-\gamma)})dq_k\wedge dq_i\\
& &\ \ \ =\sum_{i}dp_i\wedge dq_i -\frac{1}{2}\sum_{k<i}
\{(\frac{\sigma'(q_{ik}+\gamma)}{\sigma(q_{ik}+\gamma)}
+\frac{\sigma'(q_{ik}-\gamma)}{\sigma(q_{ik}-\gamma)})
-(\frac{\sigma'(q_{ik}-\gamma)}{\sigma(q_{ik}-\gamma)}
+\frac{\sigma'(q_{ik}+\gamma)}{\sigma(q_{ik}+\gamma)})\}dq_k\wedge dq_i\\
& & \ \ \ =\sum_{i}dp_i\wedge dq_i,
\end{eqnarray*}
\noindent where we have used the property that the elliptic function 
$\sigma(u)$ is an odd function with regard to argument $u$.
\hspace{2cm} $\large\bf \Box$

\vspace{0.6cm}
It is well-known that the Poisson bracket is invariant under the Poisson 
map. Hence the study of the r-matrix structure for the standard
Ruijsenaars
Lax operator $L_R(u)$ is equivalent to that of Lax operator $\stackrel{\sim}
{L}_{R}(u)$.

Following the work of Nijhoff et al \cite{Nij961}, the fundamental Poisson
bracket of 
the Lax operator $\stackrel{\sim}{L}_{R}(u)$  can be given in the 
following quadratic r-matrix form with a dynamical r-matrices  
\begin{eqnarray}
&&\{\stackrel{\sim}{L}_{R}(u)_1,\stackrel{\sim}{L}_{R}(v)_2\}
=\stackrel{\sim}{L}_{R}(u)_1\stackrel{\sim}{L}_{R}(v)_2 r^{-}_{12}(u,v)
-r^{+}_{21}(v,u)\stackrel{\sim}{L}_{R}(u)_1\stackrel{\sim}{L}_{R}(v)_2
\nonumber\\
& &~~~~~~~~~~+\stackrel{\sim}{L}_{R}(u)_1s^{+}_{12}(u,v)\stackrel{\sim}{L}
_{R}(v)_2
-\stackrel{\sim}{L}_{R}(v)_2s^{-}_{12}(u,v)\stackrel{\sim}{L}_{R}(u)_1,
\label{Qua1}
\end{eqnarray}
\noindent where
\begin{eqnarray*}
& & r^{-}_{12}(u,v)=a_{12}(u,v)-s_{12}(u)+s_{21}(v),\ \
r^{+}_{12}(u,v)=a_{12}(u,v)+u^{+}_{12}+u^{-}_{12},\\
& & s^{+}_{12}(u,v)=s_{12}(u)+u^{+}_{12},\ \ \
 s^{-}_{12}(u,v)=s_{21}(v)-u^{-}_{12},
\end{eqnarray*}
and
\begin{eqnarray*}
& & a_{12}(u,v)=r^{0}_{12}(u,v)+\sum_{i=1}\xi(u-v)e_{ii}\otimes e_{ii}
+\sum_{i\neq j}\xi(q_{ij})e_{ii}\otimes e_{jj},\\
& & r^{0}_{12}(u,v)=\sum_{i\neq j}\frac{\sigma(q_{ij}+u-v)}{\sigma(q_{ij})
\sigma(u-v)}e_{ij}\otimes e_{ji},\ \ \
s_{12}(u)=\sum_{i,j}\left( \stackrel{\sim}{L}_{R}(u)\partial_{\gamma}
\stackrel{\sim}{L}_{R}(u)\right)^{i}_{j}e_{ij}\otimes e_{jj},\\
&& u^{\pm}_{12}=\sum_{i,j}\xi (q_{ji}\pm \gamma)e_{ii}\otimes e_{jj}.
\end{eqnarray*}
The matrix element of $e_{ij}$ is equal to $(e_{ij})^{l}_{k}
=\delta_{il}\delta_{jk}$. 
It can be checked that the following symmetric condition
hold for the r-matrices $r^{\pm}_{12}(u,v)$ and $s^{\pm}_{12}(u,v)$
\begin{eqnarray}
& & r^{\pm}_{21}(v,u)=-r^{\pm}_{12}(u,v)\ \ \ ,\ \ \
s^{+}_{21}(v,u)=s^{-}_{12}(u,v),\label{Skew3}\\
& &
r^{+}_{12}(u,v)-s^{+}_{12}(u,v)=r^{-}_{12}(u,v)-s^{-}_{12}(u,v).\label{Skew4}
\end{eqnarray}
The classical r-matrices $r^{\pm}_{12}(u,v)$,
$s^{\pm}_{12}(u,v)$ are of  dynamical ones (i.e the matrix element of theirs
do depend upon the dynamical variables $q_i$). The quadratic Poisson
bracket (\ref{Qua1}) 
and the symmetric conditions of  (\ref{Skew3})-(\ref{Skew4}) lead to the
evolution integrals $tr(\stackrel{\sim}{L}_{R}(u))^{l}$.

Due to the r-matrices depending on the dynamical variables, the Poisson
bracket of $\stackrel{\sim}{L}_{R}(u)$ is no longer closed . The
complexity of the r-matrices (\ref{Qua1}) results in that it is still an
open problem to check the generalized Yang-Baxter relations for the RS
model. Moreover, the quantum version of the algebric relation 
(\ref{Qua1}) is still 
not found. The same situation also  occurs for the standard Lax operator 
$L_{R}(u)$ , and the corresponding r-matrices was given by Suris
\cite{Sur96}.

\section{The ``good" Lax representaion of elliptic RS model and its r-matrix}
The L-operator of the elliptic RS model given by Ruijsenaars $L_{R}(u)$ in
(\ref{Ru})  (or its Poisson equivalent counterpart
$\stackrel{\sim}{L}_{R}(u)$ 
in (\ref{Lr})) and corresponding r-matrix $r_{12}(u,v)$ given by Suris
\cite{Sur96} (or given by Nijhoff et al \cite{Nij961}  )
leads to some difficulties  in the investigation of the RS model. 
This motivates us to find a ``good" Lax representation of the RS model. 
As see from proposition 1 and corollary 1 in section II, this means 
to find   $g(u)$  which satisfies
(\ref{Con1})---(\ref{Con2}). 
In our former work \cite{Hou99},  
we have succeeded in find such a $g(u)$ for the elliptic RS model with 
 a special case $n=2$. 
Fortunately, we could also find such a $g(u)$ for the elliptic RS model with 
 a generic $n$ $(n>2)$ (This kind L-operator does not always exist for
general completely integrable system). The fundamental Poisson bracket of
this new L-operator  
$L(u)$ would be expressed in the Sklyanin Poisson-Lie bracket form with a 
numeric r-matrix. The corresponding r-matrix is the  classical 
$Z_n$-symmetric r-matrix  in \cite{Hou97}. Namely,  
the elliptic RS and the corresponding non-relativistic version---the 
elliptic $A_{n-1}$ CM model \cite{Hou97} are governed by the exact same
r-matrix 
(cf.\cite{Sur96}) in some gauge. In order to compare with the L-operator
given
by Ruijsenaars $L_{R}(u)$ and its Poisson equivalence
$\stackrel{\sim}{L}_{R}(u)$, we call
this L-operator  as the new Lax operator (alternatively, a ``good" 
Lax  operator).

Set an $n\otimes n$ matrix $A(u;q)$
\begin{eqnarray}
& &A(u;q)^{i}_{j}\equiv A(u;q_1,q_2,...,q_n)^{i}_{j}=
\theta^{(i)}(u+nq_j-\sum_{k=1}^{n}q_k+\frac{n-1}{2}).
\end{eqnarray}
We remark  that $A(u,q)^{i}_{j}$ correspond to
the interwiner function $\varphi^{(i)}_j$ between the $Z_n$-symmetric
Belavin model and the $A^{(1)}_{n-1}$ face model \cite{Jim88} in
\cite{Hou93}.

Define 
\begin{eqnarray*}
& &g(u)=A(u;q)\Lambda(q),\ \ \Lambda(q)^{i}_{j}=h_{i}(q)\delta^{i}_{j},
\\
& &h_{i}(q)\equiv h_{i}(q_1,....,q_n)=\frac{1}{\prod_{l\ne
i}\sigma(q_{il})}.
\end{eqnarray*}
Let us construct the new Lax operator $L(u)$ 
\begin{eqnarray}
L(u)=g(u)\stackrel{\sim}{L}_{R}(u)g^{-1}(u).
\end{eqnarray}
It will turn out that such a Lax operator $L(u)$ give a ``good" Lax
representation for the elliptic RS model. This is our main results of this
paper. To recover this, let us express the ``good" Lax operator $L(u)$
more explicitly.

\noindent {\it {\bf Proposition 3.} The Lax operator $L(u)$ can be rewritten
in the factorized form
\begin{eqnarray}
L(u)^{i}_{j}=\sum_{k=1}^{n}\frac{1}{\sigma(\gamma)}A(u+n\gamma;q)^i_k
A^{-1}(u;q)^k_je^{p_k},\ \ i,j=1,2,..., n.
\end{eqnarray}}

\noindent {\it \bf Proof:} First, let us introduce a matrix $T(u)$ with matrix
elements
\begin{eqnarray*}
T(u)^{i}_j=\sum_{k} e^{p_j}A^{-1}(u;q)^i_kA(u+n\gamma ;q)^k_j.
\end{eqnarray*}
From the definition of $A(u;q)^i_j$ and the determinant formula
of Vandermonde type \cite{Hou93}
\begin{eqnarray}
det[\theta^{(j)}(u_k)]=Const.\times \sigma (\frac{1}{n}\sum_{k}u_k
-\frac{n-1}{2})\prod_{1\leq j<k\leq n}\sigma(\frac{u_k-u_j}{n}),
\end{eqnarray}
where the Const. does not depend upon $\{u_k\}$, we have
\begin{eqnarray*}
\sum_kA^{-1}(u;q)^{i}_kA(u+n\gamma ;q)^k_j=\frac{\sigma(\gamma +u+q_{ji})}
{\sigma(u)}\prod_{k\neq i}\frac{\sigma(\gamma+q_{jk})}{\sigma(q_{ik})}.
\end{eqnarray*}
Namely,
\begin{eqnarray*}
& &T(u)^i_j=\frac{e^{p_j}\sigma(\gamma+u+q_{ji})}{\sigma(u)}
\prod_{k\neq i}\frac{\sigma(\gamma+q_{jk})}{\sigma(q_{ik})}\\
& &\ \ \ \ =\frac{1}{\prod_{k\neq i}\sigma(q_{ik})}\left\{
\frac{e^{p_j}\sigma(\gamma+u+q_{ji})\sigma(\gamma)}{\sigma(u)
\sigma(\gamma+q_{ji})}\prod_{k\neq j}\frac{\sigma(\gamma+q_{jk})}
{\sigma(q_{jk})}\right\}
\prod_{k\neq j}\sigma(q_{jk}).
\end{eqnarray*}
\noindent Then, we obtain
\begin{eqnarray*}
& &\frac{1}{\sigma(\gamma)}\sum_kA(u+n\gamma;q)^i_kA^{-1}(u;q)^k_je^{p_k}\\
& & \ \ \ =\frac{1}{\sigma(\gamma)}\sum_{m,l}A(u;q)^i_mT^m_l(u)
A^{-1}(u;q)^l_j\\
& &\ \ \  =\sum_{m,l}\frac{A(u;q)^i_m}{\prod_{k\neq m}\sigma(q_{mk})}
\left\{
\frac{e^{-p_l}\sigma(\gamma+u+q_{lm})}{\sigma(u)
\sigma(\gamma+q_{lm})}\prod_{k\neq l}\frac{\sigma(\gamma+q_{lk})}
{\sigma(q_{lk})}\right\}
A^{-1}(u;q)^l_j\prod_{k\neq l}\sigma(q_{lk})\\
& & \ \ \ =\sum_{m,l}g(u)^i_m\stackrel{\sim}{L}_{R}(u)^m_lg^{-1}(u)^l_j
\equiv L(u)^i_j
\hspace{2cm} {\bf \Box}
\end{eqnarray*}
\vspace{0.5cm}

Let us consider the non-relativistic limit of our Lax operator $L(u)$. First,
rescale the monenta $\{p_i\}$, the coupling constant $\gamma$ and the Lax
operator $L(u)$ as follows \cite{Nij961}
\begin{eqnarray}
p_i:=-\beta p'_i\ \ \ ,\ \ n\gamma :=\beta s\ \ ,\ \ L(u):=
\sigma(\frac{\beta s}{n})L'(u),
\end{eqnarray}
where $p'_i$ is the  conjugated monenta of $q_i$ in  the CM model.

Then the non-relativistic limit is obtained by taking
$\beta \longrightarrow 0$, we have the following asympotic properties
\begin{eqnarray*}
L'(u)^i_j=\delta^i_j-\beta(\sum_{k}\{A(u;q)^i_kA^{-1}(u;q)^k_jp'_k-s\partial_{u}
(A(u;q)^i_k)A^{-1}(u;q)^k_j\} ) +0(\beta^2).
\end{eqnarray*}
\noindent If we make the canonical transformation
\begin{eqnarray*}
p'_i\longrightarrow p'_i-\frac{s}{n}\frac{\partial}{\partial q_i} lnM(q),
\ \ \ M(q)=\prod_{i<j}\sigma(q_{ij}),
\end{eqnarray*}
we obtain the ``good" Lax operator of the elliptic
$A_{n-1}$ CM  model in \cite{Hou97}
\begin{eqnarray}
L_{CM}(u)^i_j=-\lim\limits_{\beta\to 0} \frac{L'(u)^i_j-\delta^i_j}{\beta}|
_{p'_i\longrightarrow p'_i-\frac{s}{n}\frac{\partial}{\partial
q_i}lnM(q)}.\label{Lcm}
\end{eqnarray}
\vspace{0.6cm}

Now, we have a position to calculate the r-matrix structure of the ``good"
Lax operator $L(u)$ for the elliptic RS model. From proposition 3 and
through the straightforward calculation, we have the main theorem of this
paper:

\noindent {\it {\bf Theorem 1.} (Main Theorem) The fundamental Poisson
bracket of $L(u)$ can be given in the quadratic Poisson-Lie form with a
nondynamical r-matrix (or Sklyanin bracket)
\begin{eqnarray}
\{L_1(u),L_2(v)\}=[r_{12}(u-v),L_1(u)L_2(v)],
\end{eqnarray}
where the numeric r-matrix $r_{12}(u)$ is the classical
$Z_{n}$-symmetric r-matrix \cite{Hou97}
\begin{eqnarray}
r^{lk}_{ij}(v)=\left\{\begin{array}{ll}
(1-\delta^{l}_{i})\frac{\theta^{'(0)}(0)\theta^{(i-j)}(v)}
{\theta^{(l-j)}(v)\theta^{(i-l)}(0)}
+\delta^{l}_{i}\delta^{k}_{j}(
\frac{\theta^{'(i-j)}(v)}{\theta^{(i-j)}(v)}-\frac{\sigma'(v)}{\sigma(v)})
&{\rm if }~i+j=l+k~{\rm mod}~ n\\
0&{\rm otherwise}
\end{array}
\right..\label{r-matrix}
\end{eqnarray}}
\vspace{0.5cm}

\noindent {\bf Remark:} I. The elliptic RS and CM model
are governed by the exact same nondynamical r-matrix in the sepcial Lax
representation.

$\ \ \ \ \ \ $ II. It was shown in \cite{Hou97} that such a
$Z_{n}$-symmetric
r-matrix satisfies the nondynamical classical Yang-Baxter equation
\begin{eqnarray}
& &[r_{12}(v_1-v_2),r_{13}(v_1-v_3)]+[r_{12}(v_1-v_2),r_{23}(v_2-v_3)]
+[r_{13}(v_1-v_3),r_{23}(v_2-v_3)]=0,\no
\end{eqnarray}
\noindent and enjoys in the antisymmetric properties
\begin{eqnarray}
-r_{21}(-v)=r_{12}(v).
\end{eqnarray}
\noindent Moreover, the r-matrix $r_{12}(u)$  also enjoys
in the $Z_{n}\otimes Z_{n}$ symmetry
\begin{eqnarray}
r_{12}(v)=(a\otimes a)r_{12}(v)(a\otimes a)^{-1}\ \ \ ,
\ \ {\rm for}\ \ a=g,h,
\end{eqnarray}
\noindent where the $n\times n$ matrices $h,g$ are defined in section 5.

\noindent {\it {\bf Corollary 2. }The Lax operator $L_{CM}(u)$
of the elliptic $A_{n-1}$ CM model in (\ref{Lcm}) 
satisfies the  nondynamcial linear Poisson-Lie bracket
\begin{eqnarray}
\{L_{CM}(u)_1,L_{CM}(v)_2\}=[r_{12}(u-v),L_{CM}(u)_1+L_{CM}(v)_2].
\label{Lin1}
\end{eqnarray}}
The direct proof that such a ``good" (classical) Lax operator $L_{CM}(u)$
of the elliptic $A_{n-1}$ CM model satisfies (\ref{Lin1})  was given 
in \cite{Hou97}.

\section{The quantum L-oprator for the elliptic quantum RS model}
In this section, we will construct the quantum L-operator for the quantum
elliptic RS model, which satisfies the nondynamical
``RLL=LLR" relation.

We first introduce the elliptic $Z_n$-symmetric quantum R-matrix related to
$Z_n$-symmetric Belavin model, which is the quantum version of the
classical $Z_n$-symmetric r-matrix defined in (\ref{r-matrix}).

We define $n\times n$ matrices $h$, $g$ and $I_{\alpha}$ by 
\begin{eqnarray*}
h_{ij}=\delta_{i+1,j{\rm mod n}},\ \ g_{ij}=\omega^{i}\delta_{i,j}
,\ \ I_{\alpha_1,\alpha_2}\equiv I_{\alpha}=g^{\alpha_2}h^{\alpha_1},
\end{eqnarray*}
where $\alpha_1,\alpha_2 \in Z_n$ and
$\omega=exp(2\pi\frac{\sqrt{-1}}{n})$.
Define the $Z_n$-symmetric Belavin's R-matrix
\cite{Ric86,Jim88,Hou93}
\begin{eqnarray}
R^{lk}_{ij}(v)=\left\{ \begin{array}{ll}
\frac{\theta^{'(0)}(0)\sigma(v)\sigma(\sqrt{-1}\hbar)}
{\sigma'(0)\theta^{(0)}(v)\sigma(v+\sqrt{-1}\hbar)}
\frac{\theta^{(0)}(v)\theta^{(i-j)}(v+\sqrt{-1}\hbar)}
{\theta^{(i-l)}(\sqrt{-1}\hbar)\theta^{(l-j)}(v)}&{\rm if }\ \
i+j=l+k\ \ {\rm mod\ \ n} \\
0&{\rm\ \ \ otherwise}
\end{array}
\right.,
\end{eqnarray}
where $\hbar$ is the Planck's constant and $\sqrt{-1}\hbar$ is 
usually called as the crossing parameter of the R-matrix. 
We  remark that our  R-matrix coincide
with the usual one in \cite{Hou93} up to a scalar factor $
\frac{\theta^{'(0)}(0)\sigma(v)}{\sigma'(0)\theta^{(0)}(v)}
\prod_{j=1}^{n-1}\frac{\theta^{(j)}(v)}{\theta^{(j)}(0)}$
, which is to make (\ref{Re}) satisfied. 
The R-matrix satisfies quantum Yang-Baxter
equation (QYBE)
\begin{eqnarray}
R_{12}(v_1-v_2)R_{13}(v_1-v_3)R_{23}(v_2-v_3)=R_{23}(v_2-v_3)R_{13}(v_1-v_3)
R_{12}(v_1-v_2).
\end{eqnarray}
Moreover, the R-matrix enjoys in following $Z_n\otimes Z_n$
symmetric properties
\begin{eqnarray}
R_{12}(v)=(a\otimes a)R_{12}(v)(a\otimes a)^{-1},\ \ {\rm for}\ \ a=g,h.
\end{eqnarray}
The $Z_n$-symmetric 
r-matrix has the following relation with its quantum counterpart
\begin{eqnarray}
& &R_{12}(v)|_{\hbar=0}=1\otimes 1,\nonumber\\
& & R_{12}(v)=1\otimes 1 +\sqrt{-1}\hbar r_{12}(v)
+0(\hbar^2)\ \ ,\ \ {\rm when\ \ }\hbar\longrightarrow 0.\label{Re}
\end{eqnarray}

Now , we construct the quantum version of L-operator $L(u)$. The usual
canonical quantization proceduce reads
\begin{eqnarray*}
p_j\longrightarrow \widehat{p}_j=-\sqrt{-1}\hbar\frac{\partial}{\partial q_j}
\ \ \ ,\ \ \ q_j\longrightarrow q_j\ \ \ ,\ \ \ j=1,....,n.
\end{eqnarray*}
Then, the corresponding quantum L-operator $\widehat{L}(u)$ consequently
reads
\begin{eqnarray}
\widehat{L}(u)^m_l
&=&\frac{1}{\sigma(\gamma)}\sum_{k=1}^{n}A(u+n\gamma;q)^m_k
A^{-1}(u;q)^k_le^{\widehat{p}_k}\nonumber\\
 &=&\frac{1}{\sigma(\gamma)}\sum_{k=1}^{n}A(u+n\gamma;q)^m_k
A^{-1}(u;q)^k_le^{-\sqrt{-1}\hbar\frac{\partial}{\partial
q_k}}.\label{Q-L}
\end{eqnarray}
\noindent It should be remarked that such a quantum L-operator is just the
factorized difference representation for the elliptic L-operator
\cite{Hou93}.
So, we have

\noindent {\it {\bf Theorem 2.} (\cite{Hou93,Qua91,Has97})
The quantum L-operator
$\widehat{L}(u)$ defined in (\ref{Q-L}) satisfies
\begin{eqnarray}
R_{12}(u-v)\widehat{L}_1(u)\widehat{L}_2(v)=\widehat{L}_2(v)\widehat{L}_1(u)
R_{12}(u-v),
\end{eqnarray}
\noindent and $R_{12}(u)$ is the $Z_n$-symmetric R-matrix.}

The proof of Theorem 2. was given by Hou et al in \cite{Hou93}, by
Quano et al in \cite{Qua91}, by Hasegawa in \cite{Has97}, through the
face-vertex
corresponding relations  independently. The direct proof was also given 
in \cite{Hou931}.

From the quantum L-operator $\widehat{L}(u)$ and the fundamental relation 
$RLL=LLR$, Hasegawa constructed the skew-symmetric fusion of 
$\widehat{L}(u)$ and succeded in  relating them with the elliptic 
type Macdonald operator in \cite{Has97}, which is actually equivalent to
the
quantum  Ruijsenaar's operators.

\section{Discussions}
In this paper, we only consider the most general RS model---the elliptic  
RS model. Such a nondynamical r-matrix structure should exist for the 
degenerated case: the rational, hyperbolic and trigonometric RS model.

From the results of the \cite{Has97,Hou971}, when the coupling
constanst
$\frac{\gamma}{\sqrt{-1}\hbar}= $ nonegative interger, the corresponding
quantum L-operator $\widehat{L}(u)$ have finite dimensional represenation.
This means that the states of quantum RS model should degenerate in
this special case.

\vskip.3in
\subsection*{ Acknowledgements.}

This work has been financially supported by National Natural Science 
Foundation of China.
We  would like to
thank Heng Fan for a careful reading of the manuscript and many
helpful comments. 
W.L.Yang was also partially supported by the grant of Northwest
University.

\vskip.3in

\section*{ Appendix. The proof of Theorem 1.}
In this appendix, we give the proof of Theorem 1, which is the main result
 of this paper.

\noindent {\it {\bf Lemma 1.} The classical L-operator $L(u)$ for the
elliptic RS model satisfies the following algebraic relations
\begin{eqnarray*}
& &[r_{12}(u-v),L_1(u)L_2(v)]^{\rho\delta}_{\alpha\beta}\\
& & \ \ =\sum_{i,j}\{A(u+n\gamma;q)^{\rho}_iA^{-1}(v;q)
^i_{\alpha}e^{p_i}\frac{\partial}{\partial
 q_i}(A(v+n\gamma;q)^{\delta}_jA^{-1}(v;q)^j_{\beta})e^{p_j}\\
& &\ \ \ \ -A(v+n\gamma;q)^{\delta}_iA^{-1}(v;q)
^i_{\beta}e^{p_i}\frac{\partial}{\partial
 q_i}(A(u+n\gamma;q)^{\rho}_jA^{-1}(u;q)^j_{\alpha})e^{p_j}\}.
\end{eqnarray*}}
\vspace{0.6cm}
\noindent {\it \bf Proof : }Let us introduce the difference operators
$\{\widehat{D}_j\}$
\begin{eqnarray*}
\widehat{D}_j=e^{-\sqrt{-1}\hbar\frac{\partial}{\partial q_j}}\ \ \ {\rm and}
\ \  \widehat{D}_j
f(q)=f(q_1,..,q_{j-1},q_j-\sqrt{-1}\hbar,q_{j+1},...,q_n).
\end{eqnarray*}
Define
\begin{eqnarray*}
T(i,j)^{\rho\delta}_{\alpha\beta}=\left\{
\begin{array}{l}
\sum_{\rho',\delta'}R(u-v)^{\rho\delta}_{\rho'\delta'}A(u+n\gamma;q)
^{\rho'}_{i}A^{-1}(u;q)^{i}_{\alpha}\widehat{D}_{i}(A(v+n\gamma;q)
^{\delta'}_{i}A^{-1}(v;q)^{i}_{\beta}), ~{\rm if }~ i=j\\
\\
\sum_{\rho',\delta'}R(u-v)^{\rho\delta}_{\rho'\delta'}\{
A(u+n\gamma;q)^{\rho'}_{i}A^{-1}(u;q)^{i}_{\alpha}
\widehat{D}_{i}(A(v+n\gamma;q)
^{\delta'}_{j}A^{-1}(v;q)^{j}_{\beta})\\
\ \ \ +A(u+n\gamma;q)
^{\rho'}_{j}A^{-1}(u;q)^{j}_{\alpha}\widehat{D}_{j}(A(v+n\gamma;q)
^{\delta'}_{i}A^{-1}(v;q)^{i}_{\beta})\},\ \ {\rm if}\ \
i\not= j
\end{array}
\right.,
\end{eqnarray*}
\noindent and 
\begin{eqnarray*}
G(i,j)^{\rho\delta}_{\alpha\beta}=\left\{
\begin{array}{l}
\sum_{\rho',\delta'}R(u-v)^{\rho'\delta'}_{\alpha\beta}
A(v+n\gamma;q)
^{\delta}_{i}A^{-1}(v;q)^{i}_{\delta'}\widehat{D}_{i}(A(u+n\gamma;q)
^{\rho}_{i}A^{-1}(u;q)^{i}_{\rho'}),~ {\rm if }~i=j\\
\\
\sum_{\rho',\delta'}R(u-v)^{\rho'\delta'}_{\alpha\beta}\{
A(v+n\gamma;q)
^{\delta}_{i}A^{-1}(v;q)^{i}_{\delta'}\widehat{D}_{i}(A(u+n\gamma;q)
^{\rho}_{j}A^{-1}(u;q)^{j}_{\rho'})\\
\ \ \ +A(v+n\gamma;q)
^{\delta}_{j}A^{-1}(v;q)^{j}_{\delta'}\widehat{D}_{j}(A(u+n\gamma;q)
^{\rho}_{i}A^{-1}(u;q)^{i}_{\rho'}), \ \ {\rm if }\ \ i\not= j
\end{array}
\right..
\end{eqnarray*}
\noindent The quantum L-operator $\widehat{L}(u)$ satisfying the
$``RLL=LLR" $ relation results in
\begin{eqnarray}
T(i,j)^{\rho\delta}_{\alpha\beta}=G(i,j)^{\rho\delta}_{\alpha\beta}.
\label{Eq}
\end{eqnarray}
Considering the asympotic properties when $\hbar\longrightarrow 0$
\begin{eqnarray*}
& & R_{12}(u)=1+\sqrt{-1}\hbar r_{12}(u)+ 0(\hbar^2),\\
& & \widehat{D}_j=1-\sqrt{-1}\hbar\frac{\partial}{\partial q_j}+
0(\hbar^2),
\end{eqnarray*}
\noindent we have 

\noindent {\bf I.} if $\ \ i=j \ \ $
\begin{eqnarray*}
& &T(i,j)^{\rho\delta}_{\alpha\beta}\equiv
T^{(0)}(i,j)^{\rho\delta}_{\alpha\beta}+\sqrt{-1}\hbar
T^{(1)}(i,j)^{\rho\delta}_{\alpha\beta}+0(\hbar^2)\\
& &=A(u+n\gamma;q)^{\rho}_{i}A^{-1}(u;q)^i_{\alpha}
A(v+n\gamma;q)^{\delta}_{i}A^{-1}(v;q)^i_{\beta}\\
& &\ \ +\sqrt{-1}\hbar\sum_{\rho',\delta'}r(u-v)^{\rho\delta}_{\rho'\delta'}
A(u+n\gamma;q)^{\rho'}_iA^{-1}(u;q)^i_{\alpha}
A(v+n\gamma;q)^{\delta'}_iA^{-1}(v;q)^i_{\beta}\\
& &\ \ -\sqrt{-1}\hbar A(u+n\gamma;q)^{\rho}_iA^{-1}(u;q)
^i_{\alpha}\frac{\partial}{\partial
 q_i}(A(v+n\gamma;q)^{\delta}_iA^{-1}(v;q)^i_{\beta})+0(\hbar^2).
\end{eqnarray*}
{\bf II.} if $\ \ i\not= j \ \ $
\begin{eqnarray*}
& &T(i,j)^{\rho\delta}_{\alpha\beta}\equiv
T^{(0)}(i,j)^{\rho\delta}_{\alpha\beta}+\sqrt{-1}\hbar
T^{(1)}(i,j)^{\rho\delta}_{\alpha\beta}+0(\hbar^2)\\
&&=A(u+n\gamma;q)^{\rho}_{i}A^{-1}(u;q)^i_{\alpha}
A(v+n\gamma;q)^{\delta}_{j}A^{-1}(v;q)^j_{\beta}\\
&&~~~~~~+A(u+n\gamma;q)^{\rho}_{j}A^{-1}(u;q)^j_{\alpha}
A(v+n\gamma;q)^{\delta}_{i}A^{-1}(v;q)^i_{\beta}\\
& &~~~~~~+\sqrt{-1}\hbar
\sum_{\rho',\delta'}r(u-v)^{\rho\delta}_{\rho'\delta'}
\{A(u+n\gamma;q)^{\rho'}_iA^{-1}(u;q)^i_{\alpha}
A(v+n\gamma;q)^{\delta'}_jA^{-1}(v;q)^j_{\beta}\\
& &~~~~~~+A(u+n\gamma;q)^{\rho'}_jA^{-1}(u;q)^j_{\alpha}
A(v+n\gamma;q)^{\delta'}_iA^{-1}(v;q)^i_{\beta}\}\\
& &~~~~~~ -\sqrt{-1}\hbar \{A(u+n\gamma;q)^{\rho}_iA^{-1}(u;q)
^i_{\alpha}\frac{\partial}{\partial
 q_i}(A(v+n\gamma;q)^{\delta}_jA^{-1}(v;q)^j_{\beta})\\
& &~~~~~~+ A(u+n\gamma;q)^{\rho}_jA^{-1}(u;q)
^j_{\alpha}\frac{\partial}{\partial
 q_j}(A(v+n\gamma;q)^{\delta}_iA^{-1}(v;q)^i_{\beta})\} +0(\hbar^2).
\end{eqnarray*}
\noindent {\bf III.} if $\ \ i=j \ \ $
\begin{eqnarray*}
& &G(i,j)^{\rho\delta}_{\alpha\beta}\equiv
G^{(0)}(i,j)^{\rho\delta}_{\alpha\beta}+\sqrt{-1}\hbar
G^{(1)}(i,j)^{\rho\delta}_{\alpha\beta}+0(\hbar^2)\\
& &=A(u+n\gamma;q)^{\rho}_{i}A^{-1}(u;q)^i_{\alpha}
A(v+n\gamma;q)^{\delta}_{i}A^{-1}(v;q)^i_{\beta}\\
& &\ \ +\sqrt{-1}\hbar\sum_{\rho',\delta'}r(u-v)^{\rho'\delta'}_{\alpha\beta}
A(v+n\gamma;q)^{\delta}_iA^{-1}(v;q)^i_{\delta'}
A(u+n\gamma;q)^{\rho}_iA^{-1}(u;q)^i_{\rho'}\\
& &\ \ -\sqrt{-1}\hbar A(v+n\gamma;q)^{\delta}_iA^{-1}(v;q)
^i_{\beta}\frac{\partial}{\partial
 q_i}(A(u+n\gamma;q)^{\rho}_iA^{-1}(u;q)^i_{\alpha})+0(\hbar^2).
\end{eqnarray*}
{\bf IV.} if $\ \ i\not= j \ \ $
\begin{eqnarray*}
& &G(i,j)^{\rho\delta}_{\alpha\beta}\equiv
G^{(0)}(i,j)^{\rho\delta}_{\alpha\beta}+\sqrt{-1}\hbar
G^{(1)}(i,j)^{\rho\delta}_{\alpha\beta}+0(\hbar^2)\\
& &=A(u+n\gamma;q)^{\rho}_{i}A^{-1}(u;q)^i_{\alpha}
A(v+n\gamma;q)^{\delta}_{j}A^{-1}(v;q)^j_{\beta}\\
&&~~~~~~+A(u+n\gamma;q)^{\rho}_{j}A^{-1}(u;q)^j_{\alpha}
A(v+n\gamma;q)^{\delta}_{i}A^{-1}(v;q)^i_{\beta}\\
& &~~~~~~+\sqrt{-1}\hbar\sum_{\rho',\delta'}
r(u-v)^{\rho'\delta'}_{\alpha\beta}
\{A(v+n\gamma;q)^{\delta}_iA^{-1}(v;q)^i_{\delta'}
A(u+n\gamma;q)^{\rho}_jA^{-1}(u;q)^j_{\rho'}\\
& &~~~~~~+A(v+n\gamma;q)^{\delta}_jA^{-1}(v;q)^j_{\delta'}
A(u+n\gamma;q)^{\rho}_iA^{-1}(u;q)^i_{\rho'}\}\\
& &~~~~~~-\sqrt{-1}\hbar \{A(v+n\gamma;q)^{\delta}_iA^{-1}(v;q)
^i_{\beta}\frac{\partial}{\partial
 q_i}(A(u+n\gamma;q)^{\rho}_jA^{-1}(u;q)^j_{\alpha})\\
& &~~~~~~+ A(v+n\gamma;q)^{\delta}_jA^{-1}(v;q)
^j_{\beta}\frac{\partial}{\partial
 q_j}(A(u+n\gamma;q)^{\rho}_iA^{-1}(u;q)^i_{\alpha})\} +0(\hbar^2).
\end{eqnarray*}
\noindent Noting (\ref{Eq}) and  considering the term of the first order
with regard to $\hbar$, we have
\begin{eqnarray}
T^{(1)}(i,j)^{\rho\delta}_{\alpha\beta}=G^{(1)}(i,j)
^{\rho\delta}_{\alpha\beta}.\label{T1}
\end{eqnarray}
\noindent Multiplying by $e^{p_i+p_j}$ from the both sider of (\ref{T1})
and sum up for $i$ and $j$, we have
\begin{eqnarray*}
\sum_{i,j}T^{(1)}(i,j)^{\rho\delta}_{\alpha\beta}e^{p_i}e^{p_j}=
\sum_{i,j}G^{(1)}(i,j)^{\rho\delta}_{\alpha\beta}e^{p_i}e^{p_j}.
\end{eqnarray*}
\noindent Due to the commutativity of $\{e^{p_j}\}$ , we obtain
\begin{eqnarray*}
& &\sum_{\rho',\delta',i,j}\{r(u-v)^{\rho\delta}_{\rho'\delta'}
\{A(u+n\gamma;q)^{\rho'}_iA^{-1}(u;q)^i_{\alpha}e^{p_i}
A(v+n\gamma;q)^{\delta'}_jA^{-1}(v;q)^j_{\beta}e^{p_j}\\
& &\ \ \ \ -r(u-v)^{\rho'\delta'}_{\alpha\beta}
A(v+n\gamma;q)^{\delta}_iA^{-1}(v;q)^i_{\delta'}e^{p_i}
A(u+n\gamma;q)^{\rho}_jA^{-1}(u;q)^j_{\rho'}e^{p_j}\}\\
& & \ \ =\sum_{i,j}\{A(u+n\gamma;q)^{\rho}_iA^{-1}(v;q)
^i_{\alpha}e^{p_i}\frac{\partial}{\partial
 q_i}(A(v+n\gamma;q)^{\delta}_jA^{-1}(v;q)^j_{\beta})e^{p_j}\\
& &\ \ \ \ -A(v+n\gamma;q)^{\delta}_iA^{-1}(v;q)
^i_{\beta}e^{p_i}\frac{\partial}{\partial
 q_i}(A(u+n\gamma;q)^{\rho}_jA^{-1}(u;q)^j_{\alpha})e^{p_j}\}.
\end{eqnarray*}
Namely, we have
\begin{eqnarray*}
& &[r_{12}(u-v),L_1(u)L_2(v)]^{\rho\delta}_{\alpha\beta}\\
& & \ \ =\sum_{i,j}\{A(u+n\gamma;q)^{\rho}_iA^{-1}(v;q)
^i_{\alpha}e^{p_i}\frac{\partial}{\partial
 q_i}(A(v+n\gamma;q)^{\delta}_jA^{-1}(v;q)^j_{\beta})e^{p_j}\\
& &\ \ \ \ -A(v+n\gamma;q)^{\delta}_iA^{-1}(v;q)
^i_{\beta}e^{p_i}\frac{\partial}{\partial
 q_i}(A(u+n\gamma;q)^{\rho}_jA^{-1}(u;q)^j_{\alpha})e^{p_j}\}.
\hspace{1cm}{\bf \Box}
\end{eqnarray*}
\vspace{0.5cm}

\noindent Now, we have a position to calculate the fundamental Poisson
bracket of $L(u)$
\begin{eqnarray*}
& &\{L_1(u),L_2(v)\}^{\rho\delta}_{\alpha\beta}=\{
L(u)^{\rho}_{\alpha},L(v)^{\delta}_{\beta}\}\\
& &\ \ =\{\sum_{i}A(u+n\gamma;q)^{\rho}_iA^{-1}(u;q)^i_{\alpha}e^{p_i},
\sum_{j}A(v+n\gamma;q)^{\delta}_jA^{-1}(v;q)^j_{\beta}e^{p_j}\}\\
& &\ \ =\sum_{i,j}
\{A(u+n\gamma;q)^{\rho}_iA^{-1}(u;q)^i_{\alpha}e^{p_i}
\frac{\partial}{\partial q_i}(A(v+n\gamma;q)^{\delta}_j
A^{-1}(v;q)^j_{\beta})e^{p_j}\}\\
& &\ \ \ -A(v+n\gamma;q)^{\delta}_iA^{-1}(v;q)^i_{\beta}e^{p_i}
\frac{\partial}{\partial q_i}(A(u+n\gamma;q)^{\rho}_j
A^{-1}(u;q)^j_{\alpha})e^{p_j}\}\\
& &\ \ =[r_{12}(u-v),L_1(u)L_2(v)]^{\rho\delta}_{\alpha\beta}.
\end{eqnarray*}
\noindent We have used the Lemma 1 in the last equation. Thus, we have
\begin{eqnarray*}
\{L_1(u),L_2(v)\}=[r_{12}(u-v),L_1(u)L_2(v)].
\end{eqnarray*}

\end{document}